\documentstyle[12pt]{article}
\begin{document}
\newtheorem{theorem}{Theorem}
\newtheorem{definition}{Definition}
\newtheorem{lemma}{Lemma}
\newtheorem{proposition}{Proposition}
\begin{center}{\Large Cauchy Kernels for some Conformally Flat Manifolds}\end{center}
\begin{center}{\Large John Ryan,\\ Department of Mathematics, University of Arkansas,\\ Fayetteville, AR 72701, USA}\end{center}
\begin{quote}
\begin{abstract}
Here we will consider examples of conformally flat manifolds that are conformally equivalent to open subsets of the sphere $S^{n}$. For such manifolds we shall introduce a Cauchy kernel and Cauchy integral formula for sections taking values in a spinor bundle and annihilated by a Dirac operator, or generalized Cauchy-Riemann operator. Basic properties of this kernel are examined, in particular we examine links to singular integral operators and Hardy spaces.
\end{abstract}
\end{quote}
\begin{center}{\Large Introduction}\end{center}

\ In a number of papers the author and his collaborators, \cite{kr1,kr2,lr1,lr2,r1,r2} use conformal transformations to develop aspects of Clifford analysis, function theory, potential theory and classical harmonic analysis over certain examples of conformally flat manifolds. As pointed out in \cite{sy} conformally flat manifolds are ones that have an atlas whose transition functions are M\"{o}bius transformations. Using Cayley transformations from $R^{n}$ into $S^{n}$ one may see that the sphere $S^{n}$ is a conformally flat manifold. Also by factoring out a simply connected subdomain of either $R^{n}$ or $S^{n}$ by a Kleinian group that acts on the domain discontinuously one may construct other examples of conformally flat manifolds. These examples include cylinders, tori, real projective spaces and $S^{1}\times S^{n-1}$. Cauchy kernels for Dirac operators defined for spinor bundles over all these examples of manifolds have been constructed in \cite{kr1,kr2,lr1,lr2,r1,r2}. See also \cite{v} for the special case of the sphere.

\ While Cauchy kernels and Cauchy Integral Formulas for Dirac operators defined for general Riemannian manifolds have been introduced in several works including \cite{c,cn,m}, in the context of conformally flat manifolds one may in the examples so far considered obtain very explicit formulas for the Cauchy kernels and also Green's kernels. 

\ In this paper we turn to look at another class of examples of conformally flat manifolds where again we can obtain explicit formulas for the Cauchy kernel and other related kernels. These manifolds are obtained using M\"{o}bius transformations to glue together $n$-spheres.  If we glue together a total of $k$ such spheres to an $n$-sphere we obtain a manifold that can be seen as a sphere with $k$ bumps, or warts, on it. In turn we might also put a finite number of bumps or warts on each wart via the same process. We can continue to place finite numbers of warts on already existing warts. Again this is done via the same gluing process. If we do not glue any of the warts to more than one sphere or wart we obtain a manifold that is homeomorphic to $S^{n}$. As such this type of manifold is simply connected. However its conformal structure is more complicated than that of the sphere. By a similar process we may also glue together several copies of $R^{n}$, or copies of spheres and euclidean spaces.

\ The purpose in this paper is to give an explicit construction of the Cauchy kernels for such manifolds. This construction makes use of the unique continuation property for solutions to the Dirac equation considered here. We also show that a number of known properties of the Cauchy kernel in Euclidean space readily carry over to the context considered here. This includes a Hardy space decomposition of $L^{p}$ spaces for strongly Lipschitz hypersurfaces lying in the manifold. This Hardy space decomposition is expressed in terms of  solutions to a Dirac equation defined on the complementary domains to the hypersurface and with appropriate continuation properties to the boundary. The techniques used here mimic existing techniques from classical harmonic analysis in Euclidean space, see for instance \cite{s}.
\newline
{\bf{Acknowledgement}} The author is grateful to David Calderbank and Michael Eastwood for drawing his attention to the type of conformally flat manifolds introduced here.

\begin{center}{\Large Preliminaries}\end{center}

\ We will consider $R^{n}$ as embedded in the real $2^{n}$ dimensional Clifford algebra $Cl_{n}$ such that under Clifford algebra multiplication $x^{2}=-\|x\|^{2}$ for each $x\in R^{n}$. So each non-zero vector $x\in R^{n}$ has a multiplicative inverse given by the Kelvin inverse of $x$, namely $x^{-1}=\frac{-x}{\|x\|^{2}}$. If $e_{1},\ldots, e_{n}$ is an orthonormal basis for $R^{n}$ then $1,e_{1},\ldots,e_{n},\ldots,e_{j_{1}}\ldots e_{j_{r}},\ldots,e_{1}\ldots e_{n}$ is a basis for $Cl_{n}$. Here $1\leq r\leq n$ and $j_{1}<\ldots j_{r}$. Regarding this basis as an orthonormal basis for $Cl_{n}$ then the norm, $\|A\|$ of a vector $A\in Cl_{n}$ can be introduced in the usual way. We will also need the anti-automorphism $\sim:Cl_{n}\rightarrow Cl_{n}:\sim(e_{j_{1}}\ldots e_{j_{r}})=e_{j_{r}}\ldots e_{j_{1}}$. As usual we write $\tilde{A}$ for $\sim A$.

\ As shown in \cite{a} and elswhere any M\"{o}bius transformation $y=\psi(x)$ over the one point compatification of $R^{n}$ can be expressed as $(ax+b)(cx+d)^{-1}$ where $a$, $b$, $c$ and $d$ are elements of $Cl_{n}$ and satisfy certain constraints. Those constraints are detailed in \cite{a,kr2,p} and elsewhere. If we regard $R^{n}$ as embedded in $R^{n+1}$ and then $Cl_{n}$ is a subalgebra of $Cl_{n+1}$. Further we may regard $R^{n+1}$ as having as orthonormal basis $e_{1},\ldots,e_{n},e_{n+1}$. In this case we have the Cayley transformation $c(x)=(e_{n+1}x+1)(x+e_{n+1})^{-1}$ which maps $R^{n}$ homeomorphically to $S^{n}\backslash\{e_{n+1}\}$.

\ The Dirac operator in euclidean space is the differential operator $D=\Sigma_{j=1}^{n}e_{j}\frac{\partial}{\partial x_{j}}$. Any $Cl_{n}$ valued, differentiable function $f$ defined on a domain $U$ in $R^{n}$ is called a left Clifford holomorphic function or left monogenic function if $Df(x)=0$ on $U$. If $g$ is also defined on $U$ and takes values in $Cl_{n}$ then if $g$ is differentiable then $g$ is called a right Clifford holomorphic function, or right monogenic function if $g(x)D=0$. Here $g(x)D=\Sigma_{j=1}^{n}\frac{\partial g(x)}{\partial x_{j}}e_{j}$. Basic properties of left and right monogenic functions are described in \cite{bds} and elsewhere. The function $G(x)=\frac{x}{\|x\|^{n}}$ is an example of a function that is both left and right Clifford holomorphic. It is shown in \cite{pq} and elsewhere that if $f(y)$ is left Clifford holomorphic in the variable $y=(ax+b)(cx+d)^{-1}$ then $J(\psi,x)f(\psi(x))$ is left Clifford holomorphic in the variable $x$ where $J(\psi,x)=\frac{\widetilde{cx+d}}{\|cx+d\|^{n}}$. Moreover, \cite{pq}, $G(\psi(x)-\psi(y))=\tilde{J}(\psi,y)^{-1}G(x-y)J(\psi,x)^{-1}$.

\ In \cite{r1} the author noted that if $y=c^{-1}(x)$ where $c^{-1}$ is the inverse of the Cayley transformation then any left Clifford holomorphic function $f(y)$ defined on a domain $U$ in $R^{n}$ is transformed to a function $J(c^{-1},x)f(c^{-1}(x))$ defined on a domain $c^{-1}(U)$ lying on $S^{n}$. In \cite{r1} it is shown that such a function is annihilated by a Dirac operator $D_{s}$ acting over $S^{n}$. An explicit expression for $D_{s}$ is determined in \cite{cm}, see also \cite{lr1,r2}. For each $x\in c^{-1}(U)$ the expression $J(c^{-1},x)f(c^{-1}(x))$ takes its values in a $2^{n}$ dimensional subspace of $Cl_{n+1}$. What we are describing is a bundle structure over $S^{n}\backslash\{e_{n+1}\}$. It is now time to introduce conformally flat manifolds.

\begin{definition}
An $n$-dimensional manifold $M$ is said to be conformally flat if there is an atlas $\cal{A}$ of $M$ such that for any pair of chart maps $(\phi_{1},U_{1})$ and $(\phi_{2},U_{2})$ the transition function $\psi_{12}=\phi_{2} \phi_{1}^{-1}$ is a M\"{o}bius transformation wherever this function is defined.
\end{definition}

\ Using Cayley transformations it is clear that $S^{n}$ is an example of a conformally flat manifold.

\ Two $n$ dimensional conformally flat manifolds $M_{1}$ and $M_{2}$ are said to be conformally equivalent if there is a diffeomorphism $\Psi:M_{1}\rightarrow M_{2}$ which with respect to the atlases $\cal{A}_{1}$ and $\cal{A}_{2}$ of $M_{1}$ and $M_{2}$ that is locally a M\"{o}bius transformation.
 
\ We shall now introduce appropriate spinor bundles over conformally flat manifolds. Given a conformally flat manifold $M$ then following \cite{kr2} one can identify a pair of points $(y,Y)\in(\phi_{2}(U_{2}),Cl_{n})$ with either $(x,J(\psi_{12},x)Y)$ or $(x,-J(\psi_{12}Y)$ where $y=\psi_{12}(x)$. If a choice of signs in the previous construction can be made for each pair of chart maps in $\cal{A}$ that is compatible then we have constructed a spinor bundle $E$ over $M$ and $M$ is regarded as a conformally flat spin manifold. The sphere $S^{n}$ is an example of such a manifold. In constructing such a bundle we have used the conformal weight functions that preserve Clifford holomorphy. It follows that for any subdomain $U$ of a conformally flat spin manifold we can introduce a section $f:U\rightarrow E$ such that locally $f$ pulls back to a left Clifford holomorphic function. Such a section is called a left Clifford holomorphic section. It also follows that we can introduce a Dirac operator $D_{M}$ that locally pulls back to the euclidean Dirac operator. Moreover $D_{M}f=0$ for any left Clifford holomorphic section. 

\ For the case of the $n$-sphere we have a Cauchy kernel $G_{s}(x,y)$ for any pair of distinct points $x$ and $y$ on $S^{n}$. Explicitly $G_{s}(x,y)=\frac{x-y}{\|x-y\|^{n}}$. See \cite{lr1,r1,r2,v} for more details. 

\begin{center}{\Large Construction of the Manifolds}\end{center}

\ We begin by considering two copies, $A_{1}$ and $A_{2}$ of the annulus $A(0,\frac{1}{r},r)=\{x\in R^{n}:\frac{1}{r}<\|x\|<r\}$. We identify each point $x\in A_{1}$ with the point $-x^{-1}\in A_{2}$. This is done via a M\"{o}bius transformation $\psi':A_{1}\rightarrow A_{2}$. We now consider two copies, $S_{1}$ and $S_{2}$, of the unit sphere $S^{n}$ lying in $R^{n+1}$ Let $R_{1}$ be a copy of $R^{n}$ containg the annulus $A_{1}$ and $R_{2}$ be a copy of $R^{n}$ containing $A_{2}$. We have Cayley transformations $c_{1}:R_{1}\rightarrow S_{1}$ and $c_{2}:R_{2}\rightarrow S_{2}$. Let $C_{1}=c_{1}(A_{1})$ and $C_{2}=c_{2}(A_{2})$. In fact both $C_{1}$ and $C_{2}$ are annuli on the spheres $S_{1}$ and $S_{2}$ respectively. Let $\overline{B}_{1}$ be the copy in $R_{1}$ of the closure of the ball in $R^{n}$ centered at the origin and of radius $\frac{1}{r}$. We may similarly define $\overline{B}_{2}$ in $R_{2}$. For $i=1,2$ let $S'_{i}=S_{i}\backslash c_{i}(\overline{B}_{i})$. We adjoin $S'_{1}$ and $S'_{2}$ by identifying points in $C_{1}$ with points in $C_{2}$ via the M\"{o}bius transformation $\psi:C_{1}\rightarrow C_{2}$ where $\psi=c_{2}\psi' c_{1}^{-1}$. In this way we have used M\"{o}bius transformations to "glue"  together the two spheres $S_{1}$ and $S_{2}$. As we have only used M\"{o}bius transformations in this "gluing" process the resulting manifold $M$ is conformally flat. We denote this manifold by $S_{1}\wedge S_{2}(r)$. The reason for the $r$ apearing here is that our construction of the manifold depends on our choice of $r$ in the outer radii of the annuli $A_{1}$ and $A_{2}$. So in fact we have constructed a whole family of conformally flat manifolds.

\ It is a simple matter to see that the manifold $S_{1}\wedge S_{2}(r)$ is diffeomorphic to the sphere $S^{n}$. As such it is simply connected. As mentioned in our introduction the conformally flat manifold $S_{1}\wedge S_{2}(r)$ can be regarded as an $n$-sphere with a bump or wart.  In fact it may fairly easily be seen using dilation and inversion that $S_{1}\wedge S_{2}(r)$ is conformally equivalent to $S^{n}$. However, from our construction it may be seen that $S_{1}\wedge S_{2}(r)$ is not embedded in $R^{n+1}$. So in this sense this class of conformally flat manifolds can be viewed as not identical to $S^{n}$. Furthermore one can vary the radii of the spheres $S_{1}$ and $S_{2}$ so that one or both are no longer the unit sphere. In this process we widen the set of  conformally flat manifolds that we may consider. Furthermore we can in turn attach another sphere $S_{3}$ to either $S_{1}$ or $S_{2}$ by the same techniques that we used here to "glue" $S_{1}$ and $S_{2}$. We would denote such a manifold as $S_{1}\wedge S_{2}\wedge S_{3}(r)$. Here we will not consider attaching $S_{3}$ to both $S_{1}$ and $S_{2}$. This particular construction leads to a manifold which is a sphere with a handle. 

\ In general we may choose a finite number of points $x_{1},\ldots,x_{k}$ lying in the sphere $S^{n}$. We now excise $k$ nonintersecting closed caps $C(x_{j})$ from $S^{n}$. Each cap $C(x_{j})$ is centered at $x_{j}$, where $j=1,\ldots, k$. By the process just described using Cayley transforms we now may attach $k$ spheres $S^{1},\ldots,S^{k}$ to this open subset of $S^{n}$. Again we obtain a conformally flat manifold that is conformally equivalent to, but not identical to, the sphere $S^{n}$. 

\ Similarly we may attach copies of $R^{n}$ to $S^{n}$. To do this consider $A(0,\frac{1}{r},\infty)=\{x\in R^{n}:\|x\|>\frac{1}{r}\}$. Now using a Cayley transform $c$ we identify $A(0,\frac{1}{r},r)\subset A(0,\frac{1}{r},\infty)$ with the spherical annulus $c(A(0,\frac{1}{r},r)\subset S^{n}\backslash c(\overline{B}(0\frac{1}{r})$. In this way we "glued" one copy of $R^{n}$ to $S^{n}$. We will denote this manifold by $R_{1}\wedge S_{2}$. 

\ By picking finitely many points $x_{1},\ldots,x_{k}$ on $S^{n}$ we may adapt the process just outlined, and attach finitely many copies of $R^{n}$ to $S^{n}$. Using Kelvin inversion and other M\"{o}bius transformations one may see that such a manifold is conformally equivalent to $S^{n}\backslash\{x_{1},\ldots,x_{k}\}$.

\ Of course one may repeat this process and "glue" copies of either $R^{n}$ or $S^{n}$ to the spheres or copies of $R^{n}$ that have already been attached to $S^{n}$, and so on. If we avoid multiply gluing a sphere or euclidean space to several parts of another sphere or euclidean space we end up with a manifold that is conformally equivalent to either $S^{n}$ or $S^{n}\backslash\{x_{1},\ldots,x_{k}\}$. For these types of conformally flat manifolds we are able to construct a Cauchy kernel.

\begin{center}{\Large Construction of the Cauchy Kernel}\end{center}

\ Here we will work simply with the case $M=S_{1}\wedge S_{2}(r)$. All the other cases are relatively straightforward extensions of this one case. We will denote the kernel by $C_{M}(x,y)$. This will be a $Cl_{n+1}$ valued function defined on $M\times M\backslash\{(x,x):x\in M\}$. When $x$ and $y$ both belong to $S'_{j}$, for $j=1,2$, then $C_{M}(x,y)$ may be identified with $\frac{x-y}{\|x-y\|^{n}}$. Here we are regarding each $S'_{j}$ as being embedded in a copy of $R^{n+1}$. Let us denote the part of $S'_{1}$ that is identified with part of $S'_{2}$ by $S_{12}$, and the part of $S'_{2}$ that is identified with part of $S'_{1}$ with $S'_{2}$ by $S_{21}$. Recall that $S_{12}$ and $S_{21}$ are annuli on spheres that are identified with each other via a M\"{o}bius transformation $\psi_{12}:S_{12}\rightarrow S_{21}$. Let $M_{1}=S'_{1}\backslash S_{12}$ and $M_{2}=S'_{2}\backslash S_{21}$. Suppose now that $x\in S'_{1}$ and $y$ belongs to the part of $S_{1}\wedge S_{2}(r)$ where $S_{12}$ and $S_{21}$ are identified. Then in this case $y$ can be identified with $y_{1}\in S_{12}$ or $y_{2}\in S_{21}$ where $\psi_{12}(y_{1})=y_{2}$. In this case the kernel $C_{M}(x,y)$ is represented by $\frac{x-y_{1}}{\|x-y_{1}\|^{n}}$, which in turn may be identified with $\frac{x-\psi_{12}^{-1}(y_{2})}{\|x-\psi_{12}^{-1}(y_{2})\|^{n}}J(\psi_{12}^{-1},y_{2})$.  Furthermore if $x$ also belongs to $S_{12}$ then $x_{2}=\psi_{12}(x)$ and $C_{M}(x,y)$ can be identified with $J(\psi_{12}^{-1},x_{2})\frac{\psi_{12}^{-1}(x_{2})-\psi_{12}^{-1}(y_{2})}{\|\psi_{12}^{-1}(x_{2})-\psi_{12}^{-1}(y_{2})\|^{n}}J(\psi_{12}^{-1},y_{2})=\frac{x_{2}-y_{2}}{\|x_{2}-y_{2}\|^{n}}$.

\ Last of all we turn to the case where $x\in S'_{1}$ and $y\in S'_{2}$. In this case we note that the M\"{o}bius transformation $\psi_{12}$ has a unique continuation to a M\"{o}bius transformation $\Psi_{12}:N_{1}\rightarrow N_{2}$ where $N_{j}=S_{j}\backslash M_{j}$. In this case $y\in N_{2}$ and $y=\Psi_{12}(y_{1})$ for some $y_{1}\in N_{1}$. In this case $C_{M}(x,y)=\frac{x-\Psi_{12}^{-1}(y)}{\|x-\Psi_{12}^{-1}(y)\|^{n}}J(\Psi_{12}^{-1},y)$. 

\ While we have not covered all possibilities of locations in $x$ and $y$ on $M$ in this construction all remaining possibilities can be constructed easily by adapting the existing constructions presented so far here.

\begin{center}{\Large Concluding Remarks}\end{center}

\ Having constructed the Cauchy kernel $C_{M}(x,y)$ for $M=S_{1}\wedge S_{2}(r)$ we readily have Cauchy's integral formula.

\begin{definition}
A $(n-1)$-dimensional hypersurface $S$ lying in $S_{1}\wedge S_{2}(r)$ is called a Lipschitz hypersurface if locally $S$ is the image under M\"{o}bius transformations of Lipschitz surfaces lying in $R^{n}$. If the Lipschitz constants for these Lipschitz surfaces lying in $R^{n}$ are bounded then $S$ is said to be strongly Lipschitz.
\end{definition}

\begin{theorem}
Suppose that $U$ is an open, connected subset of $M=S_{1}\wedge S_{2}(r)$ and $f:U\rightarrow E$ is a left Clifford holomorphic section over $U$ then for each $y\in U$ and each $(n-1)$-dimensional strongly Lipschitz hypersurface $S$ lying in $U$ and bounding a subdomain of $U$ and containing $y$, then
\[f(y)=\frac{1}{\omega_{n}}\int_{S}C_{M}(x,y)n(x)f(x)d\sigma(x)\]
where $n(x)$ is the unit vector in $TM_{x}$ that is orthogonal to $TS_{x}$ and point outwards from $S$. Here $TM$ and $TS$ are the tangent bundles of $M$ and $S$ respectively. Furthermore $\sigma$ is the Lebesgue measure on $S$.
\end{theorem}

\ Having obtained this Cauchy kernel and Cauchy integral formula standard arguments developed in \cite{mc,m} and elsewhere give us the following decomposition result.

\begin{theorem}
Suppose that $S$ is a Lipschitz hypersurface lying in $M=S_{1}\wedge S_{2}(r)$ and the complement of $S$ has two components, $S^{\pm}$. Then for $1<p<\infty$ we have that
\[L^{p}(S)=H^{p}(S^{+})\oplus H^{p}(S^{-})\]
where $L^{p}(S)$ is the space of $E$ valued sections on $S$ that are $L^{p}$ integrable, and $H^{p}(S^{\pm})$ is the Hardy $p$-space of left Clifford holomorphic sections defined on $S^{\pm}$ with  non-tangential extension to $S$ belonging to $L^{p}(S)$.  
\end{theorem} 

\ The construction we gave in the previous section of the Cauchy kernel for the manifold $S_{1}\wedge S_{2}(r)$ can readily be adapted to construct similar Cauchy kernels for all the types of manifolds constructed in this paper. It follows that analogues of Theorems 1 and 2 hold in those settings too.

\end{document}